\documentclass{amsart}
\usepackage{amssymb}

\newcommand{\ZFCa}{{\operatorname{\mathsf {ZFC}}}}
\newcommand{\CH}{\operatorname{\mathsf {CH}}}

\newcommand{\rationals}{{\mathbb Q}}
\newcommand{\rest}{{\mathord{\restriction}}}
\newcommand{\add}{\operatorname{\mathsf   {add}}}

\newcommand{\dom}{{\operatorname{\mathsf {dom}}}}

\newcommand{\N}{{\mathcal N}}
\newcommand{\M}{{\mathcal M}}

\newcommand{\lft}[2]{\mathopen\ifcase#1{}\oo\or
                        \big#2\or\Big#2\else\oo\fi} 
\newcommand{\rgt}[2]{\mathclose\ifcase#1{}\oo\or
                        \big#2\or\Big#2\else\oo\fi} 
\newcommand{\SM}{{\mathcal{SM}}}

\theoremstyle{plain}
\newtheorem{theorem}{Theorem}
\theoremstyle{plain}
\newtheorem{lemma}[theorem]{Lemma}

\newtheorem{definition}[theorem]{Definition}

\begin{document}
\title{A note on duality between measure and category}
\author{Tomek Bartoszy\'{n}ski}
\address{Department of Mathematics and Computer Science\\
Boise State University\\
Boise, Idaho 83725 U.S.A.}
\thanks{Author  partially supported by 
NSF grant DMS 95-05375 and Alexander von Humboldt Foundation}
\email{tomek@math.idbsu.edu, http://math.idbsu.edu/\char 126 tomek}

\begin{abstract}
We show that there is no Erd\"os--Sierpi\'nski mapping preserving addition.  
\end{abstract}
\maketitle

Let $\M$ and $\N$ be the ideals of meager and null subsets of
$2^\omega$.
\begin{definition}
  A bijection $F : 2^\omega \longrightarrow 2^\omega $ is called
  Erd\"os--Sierpi\'nski mapping if
$$X \in \N \iff F[X] \in \M \quad \text{ and } \quad X \in \M \iff F[X] \in \N.$$
\end{definition}
\begin{theorem}[\cite{Oxt80Mea}, \cite{BJbook}]
Assume $\CH$. There exists an Erd\"os--Sierpi\'nski mapping.
  
\end{theorem}
Since the existence of Erd\"os--Sierpi\'nski mapping implies that the
ideals $\M$ and $\N$ have the same cardinal characteristics, the
existence of such mapping cannot be proved in $\ZFCa$.

Consider the space $2^\omega $ as a topological group with addition
modulo $2$.
The following question is attributed to Ryll-Nardzewski:
\begin{quote}
  Is it consistent that  there is an Erd\"os--Sierpi\'nski mapping $F$ such that 
$$\forall x,y \in 2^\omega \ F(x+y)=F(x)+F(y).$$
\end{quote}
The motivation for this question is following (see \cite{BJbook} for
more details):
\begin{definition}
  Suppose that $X \subseteq 2^\omega $. 
We say that $X \in \mathcal  {SN}$ ($X$ has strong measure zero) if
for every set $F \in \M$, $X+F \neq 2^\omega$.

$X \in \SM$  ($X$ is strongly meager) if for every $H \in
  \N$, $X+H \neq 2^\omega $.
\end{definition}

An Erd\"os--Sierpi\'nski mapping satisfying $F(x+y)=F(x)+F(y)$ would
also map strong measure zero sets onto strongly meager sets and vice
versa. 

Consider the following statement (considered by Carlson in \cite{CarStr}):
\begin{quote}
($ \varphi $) For every set $F \in \M$ there exists a set $F' \in \M$
such that 
$$\forall x_1, x_2 \in 2^\omega \ \exists x \in 2^\omega 
\lft1((2^\omega \setminus F')+x_1\rgt1) \cup \lft1((2^\omega \setminus
F')+x_2\rgt1) \subseteq (2^\omega \setminus F)+x.$$ 
\end{quote}
Let $ \varphi^\star$ be the dual statement obtained by replacing $\M$
by $\N$.

Note that $ \varphi$ implies that $ \mathcal  {SN} $ is a 
ideal and $ \varphi^\star$ implies that $ \SM $ is an ideal (see
remarks at the end of the paper). 

\begin{theorem}[Carlson, \cite{CarStr}]
  $\ZFCa \vdash \varphi $.
\end{theorem}
\begin{proof}
For completeness we present a short proof based on the
following classical characterization of meager sets in $2^\omega $.
\begin{theorem}[\cite{BJbook}]
  A set $F \subseteq 2^\omega $ is meager if and only if there exists
  a partition of $\omega $ into intervals, $\{I_n: n \in \omega\}$ and
  a function $x_F \in 2^\omega $ such that 
$$F \subseteq \{x \in 2^\omega : \forall^\infty n \ x \rest I_n \neq
x_F \rest I_n\}.$$
\end{theorem}

Suppose that $F \subseteq 2^\omega $ is a meager set. Without loss of
generality we can assume that 
$F= \{x \in 2^\omega : \forall^\infty n \ x \rest I_n \neq
x_F \rest I_n\}$ for some partition $\{I_n: n \in \omega\}$ and real
$x_F \in 2^\omega $.

Let $J_n = I_{2n} \cup I_{2n+1}$ for every $n$.
Define
$$F'= \{x\in 2^\omega : \forall^\infty n \ x \rest J_n \neq 0\}.$$
Suppose that $x_1, x_2 \in 2^\omega $. Define
$x = x_1 \rest \bigcup_{n} I_{2n} \cup x_2 \rest \bigcup_{n}
I_{2n+1}$.
It is clear that 
$$\lft1((2^\omega \setminus F')+x_1\rgt1) \cup \lft1((2^\omega \setminus
F')+x_2\rgt1) \subseteq (2^\omega \setminus F)+x.$$
A small modification of the above argument shows that we can consider
more than two translations, countably many or even $< \add(\N)$.  
\end{proof}

\begin{theorem}
  $\ZFCa \vdash \neg\varphi^\star $.
\end{theorem}
\begin{proof}
We start with the following easy observation:
\begin{lemma}\label{first}
  Suppose that $ I \subseteq \omega $ is a finite set and $J'
  \subseteq J \subseteq 2^I$ are such that 
  \begin{enumerate}
  \item $|J'| \cdot 2^{-|I|}=1- \delta$ and $|J| \cdot 2^{-|I|}=1-
    \varepsilon$,
  \item $ \delta^2 < \varepsilon < \delta$.
  \end{enumerate}
There exist $t_1, t_2 \in 2^I$ such that 
$$ \forall s \in 2^I \ (J'+t_1) \cup (J'+t_2) \not \subseteq J+s.$$
\end{lemma}
\begin{proof}
Let
$$Z=\{(t_1,t_2,z): z \in (J'+t_1) \cup (J'+t_2)\}.$$
Check that for every $z \in 2^I$,
$$ \frac{|2^I \times 2^I \setminus (Z)_z|}{2^{2\cdot |I|}} =
  \delta^2.$$
Thus $(Z)_z \cdot 2^{-2\cdot |I|}=1- \delta^2 > 1-\varepsilon$ for all $z$.
By Fubini theorem there are $t_1, t_2$ such that 
$$\frac{|(Z)_{t_1,t_2}|}{2^{|I|}} > 1-\varepsilon. $$
In particular,
$$(Z)_{t_1,t_2} = (J'+t_1) \cup (J'+t_2) \not \subseteq J+s.$$  
\end{proof}

Fix a partition of $ \omega $ into finite sets $\{I_n: n \in \omega\}$
such that $|I_n|> 2^{n}$.
For each $n$ chose $J_n \subseteq 2^{I_n}$ such that 
$$1- \frac{1}{n^{2}}+\frac{1}{n^{5}}\geq \frac{|J_n|}{2^{|I_n|}}\geq 1- \frac{1}{n^{2}}.$$ 

Let 
$$F = \{x \in 2^\omega : \exists^\infty n \ x \rest I_n \not \in J_n\}.$$

The following lemma finishes the proof.
\begin{lemma}
  For every  null set $F' \supseteq F$ there are $x_1, x_2 \in 2^\omega
  $ such that for every $x \in 2^\omega $
$$\lft2((2^\omega \setminus F')+x_1\rgt2) \cup \lft2((2^\omega
\setminus F')+x_2\rgt2) \not \subseteq (2^\omega \setminus F)+x.$$
\end{lemma}
\begin{proof}
For a closed set $C \subseteq 2^\omega $, $ n \in \omega $ and $s \in 2^{<\omega}$ let
$$C_s =\{x \rest (|s|, \omega): s \subseteq x\} \text{ and } C \rest n
= \{x \rest n: x \in C\}.$$

Let
$$C = \{x \in 2^\omega : \forall n \ x \rest I_n \in J_n\}.$$
Without loss of generality we can assume that $C$ has positive
measure. Suppose that $F'$ is a null set. Let $C'$ be a
set of positive measure such that 
$$F' \subseteq 2^\omega \setminus (C'+\rationals),$$
where
$\rationals = \{x \in 2^\omega : \forall^\infty n  \ x(n)=0\}.$ 

We will construct two reals $x_1, x_2$ such that for every $x$
$$(C'+x_1) \cup (C'+x_2) \not \subseteq (2^\omega \setminus F)+x.$$

Define by induction an increasing  sequence $\{n_k: k \in
\omega\}$ and $x_i \rest I_{n_k}$ for $i=1,2$. For $m \neq n_k$ we put
$x_i \rest I_m
=0$. 

Suppose that
$x_1 \rest I_1 \cup I_2 \cup \dots \cup I_{n_k}$ and
$x_2 \rest I_1 \cup I_2 \cup \dots \cup I_{n_k}$ are defined. We need
to define $n_{k+1}$ and $x_i \rest I_{n_{k+1}}$ for $i=1,2$.
Use the Lebesgue density theorem to find sequences $\{r^s: s \in C' \rest I_1 \cup I_2 \cup \dots \cup I_{n_k}\}$
and $\ell > n_k$ such that 
\begin{enumerate}
\item $\dom(s{}^\frown r^s) = I_1 \cup I_2 \cup \dots \cup I_{\ell}$,
\item the set $\bigcap_{s} C'_{s {}^\frown r^s}$
  has positive measure.
\end{enumerate}
For $m \geq \ell$ let 
$$J'_m= \left\{x \rest I_m: x \in \bigcap_{s} C'_{s {}^\frown
    r^s}\right\}.$$
Note that 
$$\bigcap_{s} C'_{s {}^\frown r^s} \subseteq \{x
\in 2^\omega :
\forall m \ x \rest I_m \in J'_m\}.$$
Since the set on the left-hand side has positive measure there must be
infinitely many $m$ such that 
$$\frac{|J'_m|}{2^{|I_m|}} > 1 - \frac{1}{m},$$
since $\prod_{m}\left(1-\dfrac{1}{m}\right)=0$.
Let  $n_{k+1}$ be first such $m$ that is bigger than $\ell$.
Apply the lemma to get sequences $t^{k+1}_1, t^{k+1}_2$ such that 
$$ \forall s \in 2^{I_{n_{k+1}}} \ (J'_{n_{k+1}}+t_1^{n_{k+1}}) \cup
(J'_{n_{k+1}}+t_2^{n_{k+1}}) \not 
\subseteq J_{n_{k+1}}+s.$$ 
Define $x_1 \rest I_{n_{k+1}} = t^{k+1}_1$ and 
$x_2 \rest I_{n_{k+1}} = t^{k+1}_2$. 
This completes the definition of $x_1$ and $x_2$.

Suppose that $x \in 2^\omega $ is given. Let $s_n = x \rest I_n$.
Without loss of generality we can assume 
$$ \exists^\infty k \ (J'_{n_{k+1}}+t_1^{n_{k+1}}) \not 
\subseteq J_{n_{k+1}}+s_{n_{k+1}}.$$ 
Let $U \subseteq \omega $ be the set of $k$ satisfying the requirement
above.
We will show that  $C' +x_1 \not \subseteq (2^\omega \setminus F)+x$.
For each $k$ let $u_k \in J'_{n_{k}}+t_1^{n_{k}}$ be such that 
$u_k \in (J'_{n_{k}}+t_1^{n_{k}}) \setminus (J_{n_{k}}+s_{n_{k}})$ if
possible, i.e. if $k \in U$.

Let $v_0=r^\emptyset$ and
$$v_{k+1}= \left\{
  \begin{array}{ll}
v_k {}^\frown u_{k/2}& \text{if } k \text{ is even}\\
v_k {}^\frown r^{v_k} & \text{if } k \text{ is odd}
  \end{array} \right. .$$

Let $z = \bigcup_k v_k$. Since $[v_k] \cap (C'+x_1) \neq
\emptyset$ for all $k$, it follows that $z \in C'+x_1$. 
On the other hand  $z \not \in (2^\omega
\setminus F)+x$ since 
$$ \exists^\infty k \ z \rest I_{n_k} \not\in
J_{n_{k}}+s_{n_{k}}.$$
\end{proof}
  
\end{proof}

{\bf Remarks}

\begin{enumerate}
\item The proof shows that there is no Erd\"os--Sierpi\'nski mapping
  $F$ such that 
$$ \forall X \subseteq 2^\omega \ \forall y \in 2^\omega \ \exists z
\in 2^\omega \ F[X+y]=F[X]+z.$$
\item It is consistent that $\SM$ is an ideal (of countable sets)
  (\cite{CarStr}). Continuum Hypothesis implies that $\SM$ is not an ideal (\cite{BaSh607}).
\end{enumerate}

\end{document}